\documentclass{amsart}
\usepackage[hidelinks]{hyperref}
\usepackage{cleveref}
\usepackage{cite}
\usepackage{graphicx}

\newcommand{\ddf}[1]{\frac{\partial #1}{\partial f}}
\newcommand{\twoddf}[1]{\frac{\partial^2 #1}{\partial f^2}}
\newcommand{\ddg}[1]{\frac{\partial #1}{\partial g}}
\newcommand{\twoddg}[1]{\frac{\partial^2 #1}{\partial g^2}}
\newcommand{\ddx}[1]{\frac{\partial #1}{\partial x}}
\newcommand{\twoddx}[1]{\frac{\partial^2 #1}{\partial x^2}}

\newcommand{\twoddw}[1]{\frac{\partial^2 #1}{\partial w^2}}
\newcommand{\ddy}[1]{\frac{\partial #1}{\partial y}}
\newcommand{\twoddy}[1]{\frac{\partial^2 #1}{\partial y^2}}
\newcommand{\ddt}[1]{\frac{\partial #1}{\partial t}}

\newcommand{\dd}[2]{\frac{\partial #1}{\partial #2}}
\newcommand{\ddtwo}[2]{\frac{\partial^2 #1}{\partial #2^2}}
\newcommand{\totdd}[2]{\frac{d #1}{d #2}}

\newcommand{\real}{\mathbb{R}}
\newcommand{\posreal}{\mathbb{R}_{>0}}
\newcommand{\LL}{\mathcal{L}}
\newcommand{\pac}[1]{\mathcal{P}_{ac}^+\left({#1}\right)}
\newcommand{\giniarg}[1]{\mathcal{G}\left[#1\right]}
\newcommand{\gini}{\mathcal{G}}

\newcommand{\ysmmin}[2]{\left(#1\wedge#2\right)}

\newcommand\numberthis{\addtocounter{equation}{1}\tag{\theequation}}

\newcommand*{\defeq}{\mathrel{\vcenter{\baselineskip0.5ex \lineskiplimit0pt
                     \hbox{\scriptsize.}\hbox{\scriptsize.}}}
                     =}

\title[Dynamics of Lorenz curves]{A Dynamical Equation for the Lorenz Curve: \\{\footnotesize Dynamics of incomplete moments of probability distributions arising from Fokker-Planck equations}}

\author{David W. Cohen}
\address{Department of Mathematics, Tufts University, Medford, MA}
\email{David.Cohen@tufts.edu}
\thanks{DWC was supported by the NDSEG Fellowship during the period that this research was performed.}

\author{Merek Johnson}
\address{Department of Mathematics, Tufts University, Medford, MA}
\email{Merek.Johnson@tufts.edu}

\author{Bruce M. Boghosian}
\address{Department of Mathematics, Tufts University, Medford, MA}
\email{Bruce.Boghosian@tufts.edu}

\date{June 10, 2025}

\begin{document}

\begin{abstract}
Fokker-Planck equations (forward Kolmogorov equations) evolve probability densities in time from an initial condition. For distributions over the real line, these evolution equations can sometimes be transformed into dynamics over the incomplete zeroth and first moments. We call this perspective the Lorenz dynamics of the system after the Lorenz curve description of distributions of wealth. This offers the benefit of presenting the dynamics over a compact domain. The integral transformation is motivated and then stated for a general class of Fokker-Planck equations. Following this, the transformed equation is solved for the heat equation and some variants thereof. Finally, some equations arising from the application of kinetic theory to idealized economic systems are transformed and analyzed in this new light.
\end{abstract}


\maketitle

\section{Introduction}\label{sec:intro}

A common characterization of a distribution of wealth is the Lorenz curve, wherein the cumulative fraction of wealth is plotted as a function of the cumulative fraction of economic agents, ordered by wealth.  This parametric curve is particularly relevant in the study of economic inequality, and several inequality metrics may be understood geometrically using the Lorenz curve, including the Gini coefficient and Hoover index.  While perhaps less natural, such a parametric plot may be constructed for the zeroth and first (or higher) incomplete moments of an arbitrary one-dimensional distribution, presuming all relevant moments exist.

In recent work \cite{MR4430373,MR4804189, MR3872473}, the evolution of wealth distributions has been studied via Fokker-Planck equations derived from micro-transactional agent-based models.  Motivated by this context, in this paper we detail the transformation between probability density and Lorenz curve and derive the evolution equation for the Lorenz dynamics.  We do this in the more general setting of McKean-Vlasov type Fokker-Planck equations, whose drift and diffusion coefficients may have a functional dependency on the density in question.  Our focus is primarily on the first incomplete moment, viewed as a function of the zeroth incomplete moment, although there are possible relevant physical reasons for studying Lorenz-type dynamics of higher moments, such as the kinetic energy for a plasma kinetic equation.

In \cref{sec:genTrans}, we present the transformation and the derivation of the evolution equation for the Lorenz curve corresponding to a positive, one-dimensional density governed by a general Fokker-Planck equation. A transformation of this type, where the new independent variable has a functional dependence on the density, is unusual and the typical procedure of defining the new differentials is not available. A portion of the transformation can be viewed as an extended hodograph transformation \cite{MR1005504, MR4492787} though the totality of the transformation is more than this. A related body of work is \cite{MR2059493,MR2255281,MR2163983}, in which equations of motion are derived for the pseudo-inverse function of an evolving one-dimensional probability density. Their motivation for pursuing the dynamics of the pseudo-inverse function is to analyze convergence rates in Wasserstein metrics, which in one-dimension can be expressed using pseudo-inverses.

In \cref{sec:solvable}, we provide worked examples for the heat equation with and without an Ornstein-Uhlenbeck drift and delta-distribution initial conditions.  Furthermore, we illustrate the transformation between such equations, which can be useful in characterizing the long-term asymptotic behavior of the heat equation dynamics.  This section serves to give some intuition for the behavior of solutions in the Lorenz domain but makes clear that the transformations of \cref{sec:genTrans} likely do not yield methods for explicitly solving general one-dimensional drift-diffusion equations. Rather, the transformed dynamics are interesting in their novelty and previously unreported in the literature.

In \cref{sec:kineticModels}, we derive the evolution equation for the Lorenz dynamics of an asset-exchange model arising in the field of econophysics, and our point of origin in this work.  We also discuss the Gini coefficient, a common inequality metric, and how its evolution may be more directly understood via the Lorenz dynamics. This is the first instance, to the authors' knowledge, that the dynamics of a system that evolves a distribution of wealth have been transformed into the associated dynamics of the economically-meaningful Lorenz curve. 

\section{The transformation in general}\label{sec:genTrans}
Let $\pac{\real}$ be the space of positive, absolutely continuous Borel probability measures over $\real$. Consider a general one-dimensional Fokker-Planck equation of the form \begin{equation}
    \ddt{\rho(x,t)} = -\ddx{}\left[\Sigma[x,t,\rho(x,t)] \rho(x,t)\right] + \twoddx{}\left[D[x,t,\rho(x,t)] \rho(x,t)\right],
\label{eq:genFPE}\end{equation} where $\Sigma: \real \times \posreal \times \pac{\real} \rightarrow \real$ is the drift coefficient and $D: \real \times \posreal \times \pac{\real} \rightarrow \posreal$ is the (non-degenerate) diffusion coefficient. The initial condition $\rho_0 \in \pac{\real}$ and the boundary conditions, enforced for all $t>0$, are \[\lim_{x\rightarrow \pm \infty} \rho(x,t) = 0.\] 

The incomplete zeroth moment, i.e. the cumulative distribution function (c.d.f.), \[F(x,t) = \int_{-\infty}^x \,dy\,\rho(y,t)\] varies strictly monotonically between 0 and 1 in $x$. The incomplete first moment is \[L(x,t) = \int_{-\infty}^x \,dy\,\rho(y,t) y.\]

We focus on Fokker-Planck equations for which $F(\infty,t)$ is a conserved quantity; this corresponds to the conservation of total probability. 

At each time $t$, the Lorenz curve is the $x$-parameterized curve given by the ordered pair \[\left(F(x,t),L(x,t)\right) \in [0,1] \times \real.\] Define $f \defeq F(x,t)$, which is in the interval $[0,1]$. Moreover, since $F(x,t)$ is monotone in $x$, its inverse $F^{-1}(f,t) = x$ is well-defined.\footnote{We avoid consideration of the generalized inverse cumulative distribution function by assuming that the evolution equations described by \cref{eq:genFPE} evolve within the space $\pac{\real}.$} Therefore, the Lorenz curve can be represented as \[\left(f,L(F^{-1}(f,t),t)\right),\] where $f\in[0,1]$ and the inverse function is defined for each $t$.

Define $\LL(f,t) \defeq L(F^{-1}(f,t),t)$. The purpose of this section is to derive the equation of motion of the curve $\LL(f,t)$ in its independent variables $f$ and $t$. The transformation from $(x,t,\rho(x,t))$ to $(f,t,\LL(f,t))$ involves creating both the new independent variable $f$ and the new dependent variable $\LL$ via integral transforms of the dependent variable $\rho$. \footnote{The transformation of independent spatial variables from $x$ to $f$ is similar to an extended hodograph transformation \cite{MR1005504}.}

The conservation of total probability implies \begin{align*}
    0 &= \totdd{}{t}\int_\real\,dx\,\rho(x,t) \\
    &= \int_\real \,dx\, \ddt{\rho(x,t)} \\
    &= \int_\real \,dx\, \left\{-\ddx{}\left[\Sigma[x,t,\rho(x,t)] \rho(x,t)\right] + \twoddx{}\left[D[x,t,\rho(x,t)] \rho(x,t)\right]\right\}\\
    &= \left[-\Sigma[x,t,\rho(x,t)] \rho(x,t)+\ddx{}\left[D[x,t,\rho(x,t)] \rho(x,t)\right]\right]_{-\infty}^{\infty}.
\end{align*}

By the fundamental theorem of calculus, it is immediate that \[\ddx{F(x,t)} = \rho(x,t)\] and \[\ddx{L(x,t)} = \rho(x,t)x.\]

For notational purposes define $G(f,t)$ to be such that $F(G(f,t),t) = f$, i.e. $G(f,t)$ is the spatial inverse of $F$, the c.d.f. of $\rho(x,t)$, at time $t$. Likewise, $G(F(x,t),t) = x.$

\subsection{Transformation of independent variables}\label{ssec:IndepTrans}
The natural independent variables of $\rho, F,$ and $L$ will be $(x,t)$ whereas the independent variables for $G$ and $\LL$ are $(f,t)$. 

The forward transformation is given by $f = F(x,t)$ and the backwards via $x = G(f,t)$.

The Jacobian of the forward transformation is \begin{equation}
\ddx{f} = \ddx{F(x,t)} = \rho(x,t).
\label{eq:jacobianForwardTrans}\end{equation}
By the inverse function theorem, the Jacobian of the backwards transformation, \begin{equation}
\ddf{x} = \ddf{G(f,t)} = \left(\ddx{F(x,t)}\bigg|_{x=G(f,t)}\right)^{-1} = \left(\rho(x,t)\big|_{x=G(f,t)}\right)^{-1} = \frac{1}{\rho(G(f,t),t)}.
\label{eq:jacobianBackwardTrans}\end{equation}

Having these on hand will allow for changing integration over $x$ into integration over $f$ and vice versa. This is important as the Fokker-Planck equations under consideration include those of McKean-Vlasov type \cite{MR3443169} for which the drift $\Sigma$ and diffusion $D$ may have integral dependence on $\rho$.

The integral of an arbitrary function $Q$ with domain $\real\times\posreal$ transforms as \begin{equation}
\int_\real \,dx\, Q(x,t) = \int_0^1\,df\, \ddf{x} Q(G(f,t),t) = \int_0^1\,df\,\frac{Q(G(f,t),t)}{\rho(G(f,t),t)}
\label{eq:intTransForward}\end{equation} and for a function $R$ with domain $[0,1]\times \posreal$,\begin{equation}
\int_0^1\,df\,R(f,t) = \int_\real\,dx\,\ddx{f} R(F(x,t),t) = \int_\real \,dx\, \rho(x,t) R(F(x,t),t).
\label{eq:intTransBackward}\end{equation}

\subsection{Transformation of dependent variable}\label{ssec:DepTrans}
The Fokker-Planck equation gives the evolution of the dependent variable $\rho$. We now compute how \[\LL(f,t)\defeq L(G(f,t),t)\] evolves in its independent variables $(f,t)$.

The $f$ derivative of $\LL$ is simple to compute as \begin{align*}
    \ddf{\LL(f,t)}&= \dd{L(G(f,t),t)}{G(f,t)}\ddf{G(f,t)}\\
    &= (\rho(x,t)x)\big|_{x=G(f,t)}\frac{1}{\rho(G(f,t),t)}\\
    &= G(f,t). \numberthis \label{eq:fDerivLL}
\end{align*}

Differentiating \cref{eq:fDerivLL} again in $f$ and invoking \cref{eq:jacobianBackwardTrans} shows that \begin{equation}
\twoddf{\LL(f,t)} = \ddf{G(f,t)} = \frac{1}{\rho(G(f,t),t)}.
\label{eq:ffDerivLL}\end{equation}

Incidentally, \cref{eq:ffDerivLL} proves that the Lorenz curve associated to a positive density is convex.

Computing the time derivative of $\LL(f,t)$, \[\ddt{\LL(f,t)} = \dd{L(G(f,t),t)}{G(f,t)}\ddt{G(f,t)} + \ddt{L(G,t)},\] is more involved. First we compute $\ddt{G(f,t)}$ and then $\ddt{L(G,t)}$.

Since $F(G(f,t),t) = f$ implies \[\dd{F(G(f,t),t)}{G(f,t)}\ddt{G(f,t)} + \ddt{F(G,t)} = 0,\] then re-arranging and using the inverse function theorem \begin{align*}
    \ddt{G(f,t)} &= -\ddf{G(f,t)}\ddt{F(G,t)} \\
    &= -\frac{1}{\rho(G(f,t),t)} \int_{-\infty}^{G(f,t)}\,dy\, \ddt{\rho} \\
    &= \frac{1}{\rho(G(f,t),t)} \left[\Sigma[y,t,\rho] \rho(y,t)-\ddy{}\left[D[y,t,\rho] \rho(y,t)\right]\right]_{-\infty}^{G(f,t)}, \numberthis \label{eq:tDerivG}
\end{align*} where the second step follows from \cref{eq:ffDerivLL}.

Following a similar procedure and suppressing some variable dependence \begin{align*}
\ddt{L(x,t)}&= \int_{-\infty}^{x}\,dy\, \ddt{\rho(y,t)}y\\
&= \int_{-\infty}^x\,dy\,\ddy{}\left[-\Sigma\rho + \ddy{(D \rho)}\right]y\\
&= \left[\left(-\Sigma \rho + \ddy{(D\rho)}\right)y\right]^x_{-\infty} - \int_{-\infty}^x\,dy\,\left[-\Sigma\rho + \ddy{(D \rho)}\right]\\
&= \left[-D\rho + y\left(-\Sigma \rho + \ddy{(D\rho)}\right)\right]^x_{-\infty} + \int_{-\infty}^x\,dy\,\Sigma \rho.
\end{align*}

Some simplification is afforded when these results are put together and we let the terms evaluated at $-\infty$ vanish,
\begin{align*}
    \ddt{\LL(f,t)} &= \dd{L(G(f,t),t)}{G(f,t)}\ddt{G(f,t)} + \ddt{L(G,t)}\\
    &= (\rho(x,t)x)\big|_{x=G(f,t)} \ddt{G(f,t)} + \ddt{L(f,t)} \\
    &=- \rho(G(f,t),t)D[G(f,t),t,\rho(G(f,t),t)] + \left(\int_{-\infty}^{G(f,t)}\,dy\,\Sigma[y,t,\rho(y,t)]\rho(y,t) \right) .
\end{align*}

Making use of \cref{eq:intTransForward} through \cref{eq:ffDerivLL} yields
\begin{align*}
\ddt{\LL(f,t)} =  -\left(\twoddf{\LL(f,t)}\right)^{-1}D&\left[\ddf{\LL(f,t)},t,\left(\twoddf{\LL(f,t)}\right)^{-1}\right]\\
& \;\;\;\;\;\;\;\;\;\; + \int_0^f\,dg\,\Sigma\left[\ddg{\LL(g,t)},t,\left(\twoddg{\LL(g,t)}\right)^{-1}\right].
\end{align*}

To simplify notation, let \[\widetilde{\Sigma}[f,t,\LL] \defeq \Sigma\left[\ddf{\LL(f,t)},t,\left(\twoddf{\LL(f,t)}\right)^{-1}\right] \]
 and \[\widetilde{D}[f,t,\LL] \defeq D\left[\ddf{\LL(f,t)},t,\left(\twoddf{\LL(f,t)}\right)^{-1}\right].\] Using subscripts to denote partial derivatives, the transformed equation can be simply expressed as \begin{equation}
     \boxed{\ddt{\LL} = - \frac{\widetilde{D}[f,t,\LL]}{\LL_{ff}} + \int_0^f\,dg\,\widetilde{\Sigma}[g,t,\LL].}
 \label{eq:genTrans}\end{equation}

Note that the time evolution of the inverse c.d.f. $G(f,t)$ can be deduced starting from \cref{eq:tDerivG}. Using that $\ddx{} =\rho\ddf{}$ and \cref{eq:ffDerivLL}, $G$ evolves according to \begin{equation}
    \ddt{G(f,t)} = \Sigma\left[G,t,\left(G_f\right)^{-1}\right]-\ddf{}\left[\frac{D\left[G,t,\left(G_f\right)^{-1}\right]}{G_f}\right]\label{eq:inverseCDFEvolution}.
\end{equation} For example, the porous medium equation, \begin{equation}
    \ddt{\rho(x,t)} = \twoddx{\rho^m(x,t)}, \label{eq:porousMedium}
\end{equation} has the $(f,t,G)$ dynamics \begin{equation*}
    \ddt{G(f,t)} = - \left[\left(G_f\right)^{-m}\right]_f,
\end{equation*} which are studied in \cite{MR2059493,MR2255281,MR2163983}. The $(f,t,\LL)$ dynamics of \cref{eq:porousMedium}, \begin{equation*}
    \ddt{\LL(f,t)} = -\left(\LL_{ff}\right)^{-m},
\end{equation*} are yet to be studied.

\subsection{Inherited initial and boundary conditions}\label{ssec:InitBoundCond}
The initial condition $\LL_0(f)$ for \cref{eq:genTrans} is \[L(G(f,0),0) = \int_{-\infty}^{G(f,0)}\,dx\,\rho_0(x)x.\]

The left boundary condition for all $t>0$ is $\LL(0,t) = 0$. Whereas the right boundary condition is inherited from the time-dependent $\rho$-dynamics as \[\LL(1,t) = L(\infty,t).\]

In the case of a conserved, finite first moment, then $\LL(1,t)$ can be normalized to unity by a change in the independent variable $x$.

\section{Solvable examples}\label{sec:solvable}

\subsection{Heat equation}\label{ssec:heatEq}
\begin{figure}[ht]
\centering
\includegraphics[width=.55\textwidth]{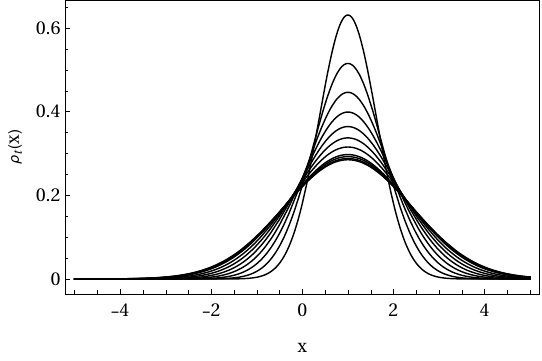}\\
\includegraphics[width=.55\textwidth]{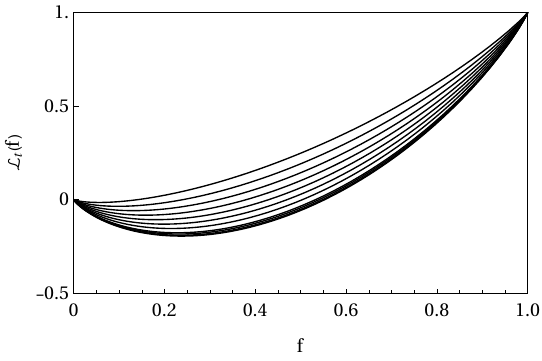}
\caption{A time-sampled solution to \cref{eq:heatEq} with $D=1$ and initial condition of $\rho_0(x) = \delta(x-1)$ (top) and the associated evolution governed by \cref{eq:heatLorenz} of its Lorenz curves with explicit solution given by \cref{eq:heatLorenzSoln} (bottom).}
\label{fig:timeVaryingHeatDistAndLor}
\end{figure}

Let $D\in\posreal$. The linear parabolic heat equation 
\begin{equation}
\ddt{\rho(x,t)} = D\twoddx{\rho(x,t)}\label{eq:heatEq}
\end{equation}
is associated to the Lorenz curve dynamics 
\begin{equation}
\LL_t = -\frac{D}{\LL_{ff}}.
\label{eq:heatLorenz}
\end{equation}
For the sake of exposition, we will first illustrate solutions for an initial condition $\rho(x,0) = \delta(x-a)$, $a\in\real$.  The fact that the first moment is conserved in the heat equation leads to Dirichlet boundary conditions, and an initial Dirac delta for $\rho$ corresponds to a linear initial condition for $\LL$,
\begin{equation}
    \begin{aligned} 
    & \LL(0,t) = 0, \quad \LL(1,t) = a, \quad  t>0,\\
    & \LL(f,0) = af, \quad 0\leq f \leq 1.
    \end{aligned}
    \label{eq:heatConditions}
\end{equation}
Solutions may be found by assuming separability of the form 
\[ 
\LL(f,t) = g(f)h(t)+k(f),
\] 
where $g$ and $h$ characterize the dynamics in the sense that $(gh)_t (gh)_{ff} = -D$ and the function $k$ serves to satisfy the initial and boundary conditions. The resulting ordinary differential equation for $h$ is easily solved, while the one for $g$ is made more tractable via the transformation $g = e^{-z^2}$ which yields \[
z' = \frac{1}{\sqrt{2}} e^{z^2}.\] This equation admits the inverse error function as a solution, and the solution to \cref{eq:heatLorenz,eq:heatConditions} is therefore
\begin{equation}
\LL(f,t) = -\sqrt{\frac{Dt}{\pi}} \exp\left(-\left(\text{erf}^{-1}(1-2f)\right)^2\right) + af.
\label{eq:heatLorenzSoln}
\end{equation}

The assumed separable form $\LL(f,t)$ is unique to the given initial condition, since it presumes that $k$ is a linear function. 

\subsection{Heat equation with Ornstein-Uhlenbeck term}
\begin{figure}[ht]
\centering
\includegraphics[width=.55\textwidth]{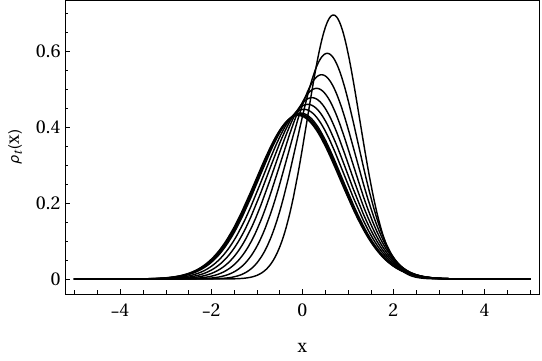}\\
\includegraphics[width=.55\textwidth]{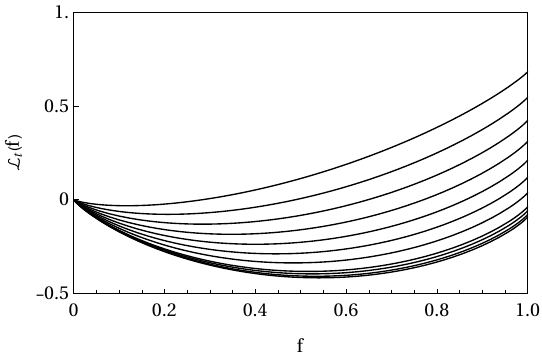}
\caption{A time-sampled solution to the \cref{eq:ouHeat} with $D,a,\sigma=1$ and $\mu=-3/4$ (top) and the associated evolution of its Lorenz curves, which obey \cref{eq:ouLorenz}, given by \cref{eq:ouLorenzSoln} (bottom). Note that the right-hand terminus of the Lorenz curve is moving downwards from the initial value of $a$ (at $t=0$) towards $\mu$ (as $t\rightarrow \infty$), in the same fashion that the peak of the Gaussian probability density shifts leftward.}
\label{fig:timeVaryingOUDistAndLor}
\end{figure}

Let $D,\sigma \in\posreal$ and $\mu\in \real$.  The heat equation with an included Ornstein-Uhlenbeck drift term is 
\begin{equation}
\ddt{\rho(x,t)} = \twoddx{}(D\rho(x,t)) +  \ddx{}\left[\sigma(\mu-x)\rho(x,t)\right]
\label{eq:ouHeat}
\end{equation}
and is associated with the Lorenz curve equation 
\begin{equation}
\LL_t = -\frac{D}{\LL_{ff}} + \sigma(\mu f - \LL). \label{eq:ouLorenz}
\end{equation}
Once again restricting to a delta initial condition for the original density centered at $a$, the left boundary and initial condition are naturally 
\begin{equation}
\begin{aligned}
    & \LL(0,t) = 0, \quad t>0\\
    & \LL(f,0) = af, \quad 0\leq f\leq 1.  
\end{aligned}
\label{eq:ouConditions}
\end{equation}
As before, the right boundary is inherited from the $\rho$-dynamics, and is now time-dependent as the mean shifts from $a$ to $\mu$. This form of the Lorenz curve equation may be solved by a transformation to the Lorenz curve equation for the heat equation, now with a nonconstant diffusion.  By defining $\mathcal{J} = \LL e^{\sigma t} - \mu f e^{\sigma t} + \alpha f +\beta$, the resulting equation for $\mathcal{J}$ becomes \[ \mathcal{J}_t\mathcal{J}_{ff} = -De^{2\sigma t}. \]
As in \cref{ssec:heatEq}, this equation is separable and the only difference manifests in the temporal ODE.  The solution to \cref{eq:ouLorenz,eq:ouConditions} is then
\begin{equation}
\LL(f,t) = -\sqrt{\frac{D(1-e^{-2\sigma t})}{2\pi\sigma}} e^{-(\text{erf}^{-1}(1-2f))^2} + (a e^{-\sigma t} + \mu (1-e^{-\sigma t})) f. 
\label{eq:ouLorenzSoln}
\end{equation}

The two analytic examples \cref{eq:heatLorenzSoln} and \cref{eq:ouLorenzSoln} of Lorenz curves are both interestingly of the form 
\begin{equation}
\LL(f,t) = -s(t) e^{-\text{erf}^{-1}(1-2f)^2} + m(t)f,
\label{eq:gensolOU}
\end{equation}
where $m(t) \defeq L(\infty,t)$ is the mean of the density $\rho$ and 
\[
s(t) \defeq \left( \int_{-\infty}^\infty dx\, \rho(x,t) (x-m(t))^2 \right)^{1/2}
\]
is the standard deviation of $\rho$.

While \cref{eq:gensolOU} seems to be a candidate for a general solution to the Lorenz dynamics of \cref{eq:ouLorenz}, it fails to match relatively simple nonlinear initial conditions such as $\LL(f,0) = f^2$. The structure of general solutions to \cref{eq:ouHeat} is given by convolving the heat kernel and the general initial condition. Despite knowledge of this convolutional structure in the $(x,t,\rho)$ domain, it has proven difficult to discern a general solution to \cref{eq:ouLorenz}.  \Cref{ssec:IndepTrans,ssec:DepTrans} demonstrate how to transform, on their own, the heat kernel and initial data; however, the correct way of combining the transformed quantities, in analogy with the convolution, in the $(f,t,\mathcal{L})$ domain is unclear and not given by \cref{eq:gensolOU}.

There are three objects that appear to demand transformation: the heat kernel, the initial condition, and the concept of a convolution. The first two are handled and the third remains an open question. Future work in this direction seems interesting.

\subsection{Transforming the quadratic potential Lorenz dynamics to the heat equation's Lorenz dynamics}\label{ssec:quadPotToHeatEq}
The asymptotic behavior of the heat equation over the real line can be difficult to understand as the time-dependent solution uniformly approaches $0$ everywhere though total heat energy is always conserved. One method of analysis is to perform a variable transformation from \[\ddt{\rho(x,t)} = \twoddx{\rho(x,t)}\] to \[\dd{\nu(y,s)}{s} =\ddy{}\left[y\nu(y,s)\right] +\twoddy{\nu(y,s)}\] and study the Gaussian asymptotic $\nu_\infty$ before transforming back to the original variables \cite[Section~2.4]{MR3497125}. The latter equation in $\nu(y,s)$ is the case of a confining quadratic potential centered on the origin.

The forward transformation is given by \[y = \frac{x}{\sqrt{2t+1}},\, s = \log\sqrt{2t+1},\, \text{ and } \nu(y,s) = e^{s}\rho\left(e^sy, \frac{1}{2}(e^{2s}-1)\right).\]

As we have seen, the $(x,t,\rho)$ system corresponds to Lorenz dynamics \[\ddt{\LL} = -\frac{1}{\LL_{ff}},\] whereas the $(y,s,\nu)$ system has Lorenz dynamics given by \begin{equation}\label{eq:quadPotLorenz}\dd{\mathcal{J}}{s} - \mathcal{J} = -\frac{1}{\mathcal{J}_{hh}}.\end{equation}

Operating wholly within the Lorenz dynamics, we derive the transformation from $(f,t,\LL)$ to $(h,s,\mathcal{J})$, which will be in the spirit of transforming $(x,t,\rho)$ to $(y,s,\nu)$.

Multiplying \cref{eq:quadPotLorenz} by $e^{2s}$ yields \[\dd{e^s\mathcal{J}(h,s)}{s} \ddtwo{e^s\mathcal{J}(h,s)}{h}=-e^{2s}\]and defining $\widetilde{\LL}(h,s) = e^s\mathcal{J}(h,s)$ reduces to \[e^{-2s}\dd{\widetilde{\LL}(h,s)}{s}\ddtwo{\widetilde{\LL}(h,s)}{h}=-1.\]

This leaves the issue of finding the temporal scaling such that $e^{-2s}\dd{}{s} = \ddt{}.$ This requires solving \[\totdd{s(t)}{t} = e^{-2s(t)}\] with \[s(t) = \log\sqrt{2t+1},\] such that $s(0) = 0.$
 
Given $s(t)$, define $\LL(h,t) = \widetilde{\LL}(h,s(t))$ in which case \[\ddt{\LL(h,t)}\ddtwo{\LL(h,t)}{h} = -1,\] yielding the transformation that we sought (trivially letting $f=h$). The transformation is \[h=f,\, s = \log\sqrt{2t+1},\, \text{ and } \mathcal{J}(h,s) = e^s\LL\left(h, \frac{1}{2}(e^{2s}-1)\right),\] which is nearly the same as the transformation in the space of densities except that now there is no need to scale the spatial variable on its own.  The reason for this is that in converting the classical heat equation to the Lorenz curve dynamics, the canonical scaling $x/2\sqrt{t}$ is already used.

\section{Kinetic models of idealized economic systems}\label{sec:kineticModels}
The motivation to pursue the transformation to Lorenz curve dynamics came from the study of idealized economic systems using the tools of kinetic theory -- a field sometimes called \textit{econophysics}, which emerged in the 1990s. In \cite{MR3428664, MR3872473}, non-linear integro-differential equations are derived that model the evolution of the distribution of wealth among a continuum of transacting economic agents. The dependent variable is $\rho(w,t),$ where the spatial variable $w>0$ is interpreted as wealth, and \[\int_a^b\,dw\,\rho(w,t)\] is the probability of choosing an agent with wealth in the interval $[a,b]$ at time $t$. 

If the time-dependent probability density $\rho(w,t)$ is such that $\int_{\posreal} \,dw\,\rho(w,t)w = 1$ then \[\int_a^b\,dw\, \rho(w,t)w\] can be interpreted as the fraction of the system's total wealth contained by agents with wealth in $[a,b]$ at time $t$.

Moreover, the systems discussed below conserve the total wealth (first moment) such that the wealth of the initial condition is the wealth for all $t>0.$ Therefore if $\int_{\posreal}\,dw\,\rho(w,0)w = 1$ then the first moment is constant and unity for all following times.

Given a stochastic transaction between two agents, deriving the governing equations for $\rho(w,t)$  under certain independence assumptions is relatively straightforward using the tools of mean-field theory or master equations.

In \cite{MR4804189}, the yard-sale model transaction is investigated: At each integer time $t$, two agents from an $N$ agent population are selected at random without replacement and their wealths are updated according to the rule \begin{equation}\label{eq:ysmMicroTransact}
\begin{pmatrix}w_{t+1}^i\\w_{t+1}^j\end{pmatrix} = \begin{pmatrix}w_{t}^i\\w_{t}^j\end{pmatrix} + \sqrt{\gamma} \ysmmin{w_t^i}{w_t^j} \begin{pmatrix}1\\-1\end{pmatrix}\eta,
\end{equation} where $\gamma\in(0,1)$ is a transaction intensity parameter, $\eta$ is a random variable with outcomes $-1$ and $+1$ with equal probability, and $\wedge$ is the $\min$ operator. If all the agents are initialized with positive wealth, then each agent will always have positive wealth.

The equation of motion for the probability distribution of agents in wealth-space under the yard-sale model dynamics of \cref{eq:ysmMicroTransact} is \begin{equation}\label{eq:ysm}\ddt{\rho(w,t)}= \twoddw{}\left[\left(\frac{\gamma}{2}\int_0^\infty\,dx\,\ysmmin{w}{x}^2\rho(x,t)\right)\rho(w,t)\right].\end{equation} We take the dynamics of \cref{eq:ysm} to occur over $\posreal\times\posreal$ when the initial condition $\rho_0(w) \in \pac{\posreal}.$ Under the continuum equation, \cref{eq:ysm}, both total probability and wealth are conserved quantities of the dynamics thus we canonically set both to be unity starting from the initial condition, \[\int_0^\infty\,dw\,\rho_0(w)w=1= \int_0^\infty\,dw\,\rho_0(w).\] We operate under these -- positive wealth and unit total wealth -- assumptions through the remainder of this section.

Thus, the boundary conditions are $\rho(0,t) = 0$ and $\lim_{w\rightarrow \infty}\rho(w,t) = 0$ for all $t>0$. This corresponds to $L(0,t) = 0$ and $L(\infty, t) = 1$ for all $t\geq 0$.

A fair portion of the study of wealth inequality investigates the evolution of the Lorenz curve rather than the distribution of agents in wealth-space. Therefore, the transformation detailed in \cref{sec:genTrans} permits a more natural study of the inequality dynamics. \footnote{Interestingly, there is no apparent (to us) derivation of the dynamics obeyed by $\LL(f,t)$ from the model's first principles, that is without passing through the derivation of the $\rho(w,t)$ dynamics.}

In \cref{eq:ysm}, \[\Sigma \equiv 0\] and \[D[w,t,\rho] = \frac{\gamma}{2}\int_0^\infty\,dx\,\ysmmin{w}{x}^2\rho(x,t).\]

Using \cref{eq:intTransForward} and \cref{eq:fDerivLL}, the diffusion coefficient transforms as \[\widetilde{D}[f,t,\LL] = \frac{\gamma}{2}\int_0^1\,dg\,\ysmmin{\LL_g(g,t)}{\LL_f(f,t)}^2.\]

By \cref{eq:genTrans}, the Lorenz curve dynamics are \begin{equation}
    \label{eq:LorenzDynYSM} \ddt{\LL} = -\frac{1}{\LL_{ff}}\frac{\gamma}{2}\int_0^1\,dg\,\ysmmin{\LL_g}{\LL_f}^2.
\end{equation}

To the best of our knowledge, a result of the type \cref{eq:LorenzDynYSM} is new to the field of econophysics. That is, we cannot find other work that takes a specified collisional transaction between two idealized economic agents and produces the dynamics of the Lorenz curve in the mean-field limit. In the context of idealized economic systems, this procedure promises a novel method to directly study the evolution of inequality.

\subsection{Intuition for the transformation via the Gini coefficient}
The Gini coefficient is a frequently used statistical measure of wealth inequality \cite{CG1912,MR3012052}. The Gini coefficient \begin{equation}\giniarg{\rho} \defeq 1- \int_0^\infty\,dw\,\int_0^\infty\,dy\,(w\wedge y)\rho(w)\rho(y)\label{eq:giniOverDist}\end{equation} is a quadratic functional of $\rho(w,t)$, when defined over the economic agent density in wealth-space.\footnote{A word of warning: This formula depends on both $\int_{\posreal}\,dw\,\rho(w) = 1 = \int_{\posreal}\,dw\,\rho(w)w$.} From the point of view of Lorenz curves of wealth distributions, the Gini coefficient has a simple geometric interpretation. Namely it is twice the area between the identity and the Lorenz curve, \begin{equation}\widetilde{\gini}\left[\LL(f,t)\right] = 2\int_0^1\,df\, \left(f-\LL(f,t)\right), \label{eq:giniOverLorenz}\end{equation} where $\widetilde{\gini}$ indicates that the functional takes an argument of a Lorenz curve whereas $\gini$ acts on the subspace of $\pac{\posreal}$ with unit first moment. For a wealth distribution $\rho$ and its Lorenz curve $\LL$ then $\giniarg{\rho} = \widetilde{\gini}\left[\LL\right].$

Let $\widetilde{\gini}(t) = \widetilde{\gini}\left[\LL(f,t)\right]$ then as a simple consequence of \cref{eq:giniOverLorenz},\begin{equation}\label{eq:dGinidt1}\totdd{\widetilde\gini(t)}{t} = -2\totdd{}{t}\int_0^1\,df\,\LL(f,t).\end{equation}

The paper \cite[Lemma~3.1]{MR4804189} proved that under the equation of motion \begin{equation}\label{eq:genEconoModel}\ddt{\rho(w,t)}= \twoddw{}\left[D[w,t,\rho]\rho(w,t)\right],\end{equation} the rate of change of the Gini coefficient $\gini(t) = \giniarg{\rho(w,t)}$ is \begin{equation}\totdd{\gini(t)}{t} = 2\int_0^\infty\,dw\,D[w,t,\rho]\rho(w,t)^2\label{eq:giniEvoUnderEcono}.\end{equation}

Transforming \cref{eq:giniEvoUnderEcono} from its expression involving an integration over $w$ into one over $f$ via \cref{eq:intTransForward} gives \begin{equation}\label{eq:dGinidt2}\totdd{\gini(t)}{t} = 2\int_0^1\,df\, \frac{\widetilde{D}}{\LL_{ff}}.\end{equation}

Let $\LL(f,t)$ be the time-evolving Lorenz curve associated to the $\rho$-dynamics generated by \cref{eq:genEconoModel}. Equating \cref{eq:dGinidt1} and \cref{eq:dGinidt2} gives \begin{equation}\totdd{}{t}\int_0^1\,df\,\LL(f,t) = -\int_0^1\,df\, \frac{\widetilde{D}}{\LL_{ff}} \label{eq:transIdenIntegrated},\end{equation} which is an integrated version of the $\Sigma \equiv 0$ case of \cref{eq:genTrans}. Of course \cref{eq:transIdenIntegrated} does not imply the stronger, pointwise statement of \cref{eq:genTrans} rather this more easily obtained result hints at the involved work in \cref{sec:genTrans}.

It has been proven in both the discrete \cite{CB2023} and continuous time \cite{MR3428664} versions of the yard-sale model that: (1) wealth condenses into a vanishingly small portion of the population and (2) the Gini coefficient approaches its maximum value of 1 in the limit $t\to\infty$. This has also been studied in a more general context in \cite{MR4267576}. For the equation of motion \cref{eq:ysm}, this increase is monotonic, a result which is easily understood in the context of \cref{eq:dGinidt2} where the convexity of the Lorenz curve is guaranteed by the positivity of the agent density $\rho$.

A direct proof of the preservation of convexity for more general Lorenz dynamics without an appeal to the behavior of $\rho$ and a minimum principle for a broader class of parabolic equations is a worthwhile problem for future study.

\section{Summary}\label{sec:summary}
We have derived a non-linear equation of motion for the partial first moment of a probability density evolving under a one-dimensional Fokker-Planck equation, possibly of the McKean-Vlasov type. The authors' original motivation to pursue this transformation came from the analysis of idealized economic systems because the time-dependent, wealth-parameterized plots of the Lorenz curve have widespread economic significance.

The Fokker-Planck equation, \cref{eq:genFPE}, with independent variables $(x,t)$ and dependent variable $\rho(x,t)$ was transformed into a system with natural variables $(f,t),$ where $f$ is associated to the cumulative distribution function of $\rho$, and dependent function $\LL(f,t)$, where $\LL(f,t)$ is the amount of the first moment contained by the fraction $f$ of the probability mass that contains the least amount of the first moment.\footnote{While this statement may be unwieldy, thinking of the economic case makes it more readily understandable: $\LL(0.5,T)$ is the fraction of wealth held by the poorest $50\%$ of a population.}

We end by mentioning two worthwhile avenues for future study: (1) Investigating how minimum and maximum principles from parabolic equations \cite{MR2597943} transform into the nonlinear Lorenz dynamics equations, and (2) Discerning how the structure of known Green's functions and their convolution in the $(x,t,\rho)$ dynamics transforms into the $(f,t,\LL)$ Lorenz dynamics. This would permit the study of more general initial conditions.

\bibliographystyle{amsplain}
\bibliography{biblio}

\end{document}